\documentclass[a4paper,leqno]{amsart}

\usepackage[utf8]{inputenc}
\usepackage[T1]{fontenc}
\usepackage{amsmath, amsthm}
\usepackage{amssymb}
\usepackage{geometry}
\usepackage{bbm}

\usepackage[hidelinks]{hyperref} 
   \hypersetup{final} 

\newtheorem{lem}{Lemma}[section]
\newtheorem{prop}[lem]{Proposition}
\newtheorem{cor}[lem]{Corollary}
\newtheorem{thm}[lem]{Theorem}

\theoremstyle{remark}
\newtheorem{rem}[lem]{Remark}

\theoremstyle{definition}

\newtheorem{example}[lem]{Example}

\numberwithin{equation}{section}

\usepackage[utf8]{inputenc}
\usepackage[T1]{fontenc}
\usepackage{amsmath, amsthm}
\usepackage{amssymb}
\usepackage{geometry}
\usepackage{color}
\usepackage{bbm}

\usepackage[hidelinks]{hyperref} 
   \hypersetup{final} 

\usepackage{mathtools}
\mathtoolsset{centercolon} 

\usepackage{todonotes} 	
\usepackage[normalem]{ulem}
\usepackage{cancel}

\def\er{{\mathbb R}}

\def\Ex{{\mathbb E}}
\def\Pr{\mathbb P}

\def\ve{\varepsilon}
\def\ind{\mathbbm{1}}

\def\Log{\mathrm{Log}}

\def\calI{\mathcal{I}}

\begin{document}

\title{On the spectral norm of Rademacher matrices}
\author{Rafał Latała}
\address{Institute of Mathematics, University of Warsaw, Banacha 2, 02--097 Warsaw, Poland.}
\email{rlatala@mimuw.edu.pl}


\subjclass[2020]{Primary 60B20, Secondary 15B52, 46B09}

\begin{abstract}
We discuss two-sided non-asymptotic bounds for the mean spectral norm of nonhomogenous 
weighted Rademacher matrices. 
We show that the recently formulated conjecture holds up to 
$\log \log \log n$  factor for arbitrary
$n\times n$ Rademacher matrices and the triple logarithm may be eliminated for matrices
with $\{0,1\}$-coefficients. 
\end{abstract}

\maketitle

\section{Introduction and main results}

One of the basic issues of the random matrix theory are bounds on the spectral norm (largest singular value) of various families of random  matrices. 
This question is very well understood for classical ensembles of random matrices
\cite{AGZbook}, when one may use methods based on the large degree of symmetry. Recently,
a substantial progress was attained in the understanding of inhomogeneous models 
\cite{vH2017_survey}, especially in the Gaussian case \cite{LvHY,BBvH}. However,
there are still many open questions in this area, the one concerning Rademacher matrices
is discussed here.

In this paper we investigate the mean operator (spectral) norm of weighted Rademacher matrices, i.e., quantities of the form
\[
\Ex\|(a_{i,j}\ve_{i,j})\|
:= \Ex\sup_{\|s\|_2,\|t\|_2 \leq 1}\sum_{i,j}a_{i,j}\ve_{i,j}s_it_j,
\] 
where $(a_{i,j})$ is a deterministic matrix and $(\ve_{i,j})_{i,j\geq 1}$ is the double indexed sequence of i.i.d.\ symmetric $\pm 1$ r.v's.

Since operator norm is bigger than length of every column and row we get
\[
\Ex\|(a_{i,j}\ve_{i,j})\|\sim (\Ex\|(a_{i,j}\ve_{i,j})\|^2)^{1/2}
\geq \max\Bigl\{\max_{i}\|(a_{i,j})_j\|_2,\max_j\|(a_{i,j})_i\|_2\Bigr\}.
\]
For two nonnegative functions 
 $f$ and $g$ we write $f \gtrsim$ g (or $g \lesssim f$) if there exists an absolute
constant $C$ such that $Cf \geq  g$; the notation $f\sim g$ means that $f \gtrsim g$ 
and $g \gtrsim f$.
Seginer \cite{Seginer} proved that for $n\geq 2$,
\[
\Ex\|(a_{i,j}\ve_{i,j})_{i,j\leq n}\|\lesssim \log^{1/4}n 
\Bigr(\max_{i}\|(a_{i,j})_j\|_2+\max_j\|(a_{i,j})_i\|_2\Bigl)
\]
and constructed an example showing that in general the constant $\log^{1/4} n$ cannot be improved. 

Let $g_{i,j}$ be independent $\mathcal{N}(0,1)$ random variables. We have
\begin{equation}
\label{eq:BvH0}
\Ex\|(a_{i,j}\ve_{i,j})\|\lesssim \Ex\|(a_{i,j}g_{i,j})\|\lesssim \max_{i}\|(a_{i,j})_j\|_2+\max_j\|(a_{i,j})_i\|_2
+\sqrt{\log n}\max_{i,j}|a_{i,j}|,
\end{equation}
where the last bound was showned by Bandeira and van Handel \cite{BvH}. 
The bound on the Gaussian
matrices is sharp for $\{0,1\}$-weights; however, even in this case, estimate \eqref{eq:BvH0} often 
yields suboptimal bounds (see Section \ref{sec:examples} below).

In \cite[Theorem 1.1]{LS} it was shown that for any  matrix $(a_{ij})$,
\begin{align}
\label{eq:loweropnorm}
\Ex\|(a_{i,j}\ve_{i,j})_{i,j\leq n}\|
&\gtrsim
\max_{1\leq i\leq n}\|(a_{i,j})_j\|_2
+\max_{1\leq j\leq n}\|(a_{i,j})_i\|_2
\\
\notag
&\phantom{aa}
+\max_{1\leq k\leq n}\min_{I\subset [n],|I|\leq k}
\sup_{\|s\|_2,\|t\|_2\leq 1}\left\|\sum_{i,j\notin I}a_{i,j}\ve_{i,j}s_it_j\right\|_{\Log\, k}.
\end{align}
Here and in the sequel $\Log\ x=\log(x\vee e)$ and
 $\|S\|_p=(\Ex|S|^p)^{1/p}$ denotes $L_p$-norm of a r.v.~$S$.

It was also conjectured that bound \eqref{eq:loweropnorm} may be reversed,
i.e., for any scalar matrix $(a_{i,j})_{i,j\leq n}$,
\begin{align}
\label{eq:upperopnorm}
\Ex\|(a_{i,j}\ve_{i,j})_{i,j\leq n}\|
&\lesssim
\max_{1\leq i\leq n}\|(a_{i,j})_j\|_2
+\max_{1\leq j\leq n}\|(a_{i,j})_i\|_2
\\
\notag
&\phantom{aa}
+\max_{1\leq k\leq n}\min_{I\subset [n],|I|\leq k}
\sup_{\|s\|_2,\|t\|_2\leq 1}\left\|\sum_{i,j\notin I}a_{i,j}\ve_{i,j}s_it_j\right\|_{\Log\, k}.
\end{align}

The proof of  \cite[Remark 4.5]{LS}, based on the permutation method from \cite{LvHY}, shows that in order to establish \eqref{eq:upperopnorm}
it is enough to show that for any submatrix 
$(b_{i,j})_{i,j\leq m}$ of $(a_{i,j})_{i,j\leq n}$ one has
\begin{equation}
\label{eq:upperopnormred}
\Ex\|(b_{i,j}\ve_{i,j})_{i,j\leq m}\|
\lesssim
\max_{1\leq i\leq m}\|(b_{i,j})_j\|_2
+\max_{1\leq j\leq m}\|(b_{i,j})_i\|_2+R_B(\Log\, m),
\end{equation} 
where for a matrix $A=(a_{i,j})$ and $p\geq 1$ we put
\[
R_A(p)
:= \sup_{\|s\|_2\leq 1,\|t\|_2\leq 1}\Big\|\sum_{i,j}a_{i,j}\ve_{i,j}s_it_j\Bigr\|_p.
\]

Our first result states that this conjectured bounds holds for $\{0,1\}$-matrices.

\begin{thm}
\label{thm:2sided01}
Inequality \eqref{eq:upperopnormred} holds if $b_{i,j}\in\{0,1\}$ for any $i,j$.
As a consequence, for any $E\subset[n]\times [n]$,
\begin{align}
\label{eq:upperopnorm01}
\Ex\|(\ind_E(i,j)\ve_{i,j})_{i,j\leq n}\|
&\sim
\max_{1\leq i\leq n}\|(\ind_E(i,j))_j\|_2
+\max_{1\leq j\leq n}\|(\ind_E(i,j))_i\|_2
\\
\notag
&\phantom{aa}
+\max_{1\leq k\leq n}\min_{I\subset [n],|I|\leq k}
\sup_{\|s\|_2,\|t\|_2\leq 1}\left\|\sum_{i,j\notin I}\ind_E(i,j)\ve_{i,j}s_it_j\right\|_{\Log\, k}.
\end{align}
\end{thm}

Inequality \eqref{eq:upperopnormred} for $\{0,1\}$-weights is a consequence of the more general Theorem 
\ref{thm:normgraphA2} below, applied to the symmetric $2m\times 2m$ \{0,1\}-matrix
$A=\left(\begin{array}{cc}0&A\\A^T&0\end{array}\right)$. Estimate \eqref{eq:upperopnorm01}
follows from \eqref{eq:upperopnormred} as in the proof of \cite[Remark 4.5]{LS}.

\begin{rem}
\label{rem:symmetrization}
In our results we do not assume the symmetry of $(\ve_{i,j})_{i,j}$ even if we consider
symmetric matrices $(a_{i,j})$.
 However analogous upper bounds holds for 
$\Ex\|(a_{i,j}\tilde{\ve}_{i,j})_{i,j}\|$, where $(\tilde{\ve}_{i,j})_{i,j}$ is the symmetric
Rademacher matrix (i.e., $\tilde{\ve}_{i,j}=\tilde{\ve}_{j,i}=\ve_{i,j}$ for $i\geq j$), since
\begin{align*}
\Ex\|(a_{i,j}\tilde{\ve}_{i,j})_{i,j}\|
&\leq \Ex\|(a_{i,j}\tilde{\ve}_{i,j}\ind_{\{i\leq j\}})_{i,j}\|
+\Ex\|(a_{i,j}\tilde{\ve}_{i,j}\ind_{\{i> j\}})_{i,j}\|
\\
&= \Ex\|(a_{i,j}\ve_{i,j}\ind_{\{i\leq j\}})_{i,j}\|
+\Ex\|(a_{i,j}\ve_{i,j}\ind_{\{i> j\}})_{i,j}\|
\leq 2\Ex\|(a_{i,j}\ve_{i,j})_{i,j}\|.
\end{align*}
Moreover, the lower bound (1.2) remains valid if we replace $\ve_{i,j}$ with $\tilde{\ve}_{i,j}$.
\end{rem}

Quantity $R_A(p)$ involves random variables $\ve_{i,j}$. It may be expressed
in terms of $a_{i,j}$ using two-sided bounds for $L_p$-norms of 
Rademacher sums (derived in \cite{Hi} on the base of tail bounds \cite{MS})
\[
\Bigg\|\sum_{k=1}^n a_k\ve_k\Bigg\|_p\sim \sum_{k\leq p}a_k^*+\sqrt{p}\Bigg(\sum_{k>p}(a_k^*)^2\Bigg)^{1/2}\sim
\sup\Bigg\{\sum_{k=1}^n a_kb_k\colon \|b\|_\infty\leq 1, \|b\|_2\leq \sqrt{p}\Bigg\},
\]
where $(a_k^*)_{k=1}^n$ denotes the nonincreasing rearrangement of $(|a_k|)_{k=1}^n$. 
It is still unclear how to apply 
the above bounds to get a simple two-sided  estimate for $R_A(p)$. However, in
\cite[Proposition 1.4]{LS}, such an estimate was established for $\{0,1\}$-matrices:
\[
\sup_{\|s\|_2,\|t\|_2\leq 1}\left\|\sum_{i,j}\ind_E(i,j)\ve_{ij}s_it_j\right\|_{p}
\sim \max_{F\subset E,|F|\leq p}\|(\ind_{\{(i,j)\in F\}})\|.
\]
Hence Theorem \ref{thm:2sided01} implies the following corollary
-- its first part provides a positive answer to the question posed by Ramon van Handel (private communication).

\begin{cor}
\label{cor:2sided}
For any $E\subset[n]\times [n]$,
\begin{align}
\label{eq:upper01}
\Ex\|(\ind_E(i,j)\ve_{i,j})_{i,j\leq n}\|
&\lesssim 
\max_{1\leq i\leq n}\|(\ind_E(i,j))_j\|_2
+\max_{1\leq j\leq n}\|(\ind_E(i,j))_i\|_2
\\
\notag
&\phantom{aa}
+\sup_{F\subset E,|F|\leq \Log\, n}\|(\ind_{\{(i,j)\in F\}})_{i,j}\|.
\end{align}
and
\begin{align}
\label{eq:2sided01}
\Ex\|(\ind_E(i,j)\ve_{i,j})_{i,j\leq n}\|
&\sim
\max_{1\leq i\leq n}\|(\ind_E(i,j))_j\|_2
+\max_{1\leq j\leq n}\|(\ind_E(i,j))_i\|_2
\\
\notag
&\phantom{aa}
+\max_{1\leq k\leq n}\min_{I\subset [n],|I|\leq k}
\max_{F\subset E,|F|\leq \Log\, k}\|(\ind_{\{(i,j)\in F,i,j\notin I\}})_{i,j}\|.
\end{align}
\end{cor}

Example \ref{ex:minoverkneeded} below shows that one cannot reverse estimate \eqref{eq:upper01}
and a more involved form \eqref{eq:2sided01} of the two-sided estimate is necessary.

\begin{rem}
Two-sided bounds on moments of norms of Rademacher vectors \cite{DMS} give that for 
every $p\geq 1$,
\[
\Bigl(\Ex\|(a_{i,j}\ve_{i,j})_{i,j\leq n}\|^p\Bigr)^{1/p}
\sim \Ex\|(a_{i,j}\ve_{i,j})_{i,j\leq n}\|+R_A(p).
\]
Thus, estimate \eqref{eq:upper01}
might be equivalently stated as
\begin{align}
\label{eq:estlognmom}
(\Ex\|(\ind_E(i,j)\ve_{i,j})_{i,j\leq n}\|^{2\lfloor \Log\ n\rfloor})
^{1/2\lfloor \Log n\rfloor}
&\sim 
\max_{1\leq i\leq n}\|(\ind_E(i,j))_j\|_2
+\max_{1\leq j\leq n}\|(\ind_E(i,j))_i\|_2
\\
\notag
&\phantom{aa}
+\max_{F\subset E,|F|\leq \Log\, n}\|(\ind_{\{(i,j)\in F\}})_{i,j}\|.
\end{align}
It is quite tempting to show \eqref{eq:estlognmom} for symmetric sets $E$ via a combinatorial
method, since for $n\times n$ symmetric matrix $A$ and 
$k=\lfloor \Log\, n\rfloor$, $\|A\|\sim (\mathrm{tr}(A^{2k}))^{1/{2k}}$. Such an approach
worked for Gaussian matrices \cite{BvH}, but we were not able to apply it in the 
Rademacher case.
\end{rem}

\begin{rem}
Signed adjacency matrices were studied in \cite{BiL} in connection with 2-lifts of graphs.
\cite[Lemma 3.1]{BiL} shows that to each signed adjacency matrix of a  graph $G$ one may associate
the 2-lift of $G$ with the set of eigenvalues being the union of the eigenvalues of $G$ and of the signed
matrix. Hence Theorem \ref{thm:2sided01} provides an average uniform bound on new eigenvalues of 
random 2-lifts.  
\end{rem}

To state results for general matrices we need to introduce some additional notation. We associate 
to a symmetric matrix $(a_{i,j})_{i,j\leq n}$
a graph $G_A=([n],E_A)$, where $(i,j)\in E_A$ iff $i\neq j$ and $a_{i,j}\neq 0$.
By $d_A$ we denote the maximal degree of vertices in $G_A$. Observe that in the case
of $\{0,1\}$-matrices 
$\sqrt{d_A}\|(a_{i,j})\|_\infty=\sqrt{d_A}= \max_{i}\|(a_{i,j})_{j}\|_2$.

\begin{thm}
\label{thm:normgraphA2}
For any symmetric matrix $(a_{i,j})_{i,j\leq n}$, 
\begin{equation}
\label{eq:graphopnormA1}
\Ex\|(a_{i,j}\ve_{i,j})_{i,j\leq n}\|
\lesssim \max_{i}\|(a_{i,j})_j\|_2+R_A(\Log\, n)+d_A^{19/40} \|(a_{i,j})\|_\infty.
\end{equation}
\end{thm}

\begin{rem}
\label{rem:zerodiag}
Since $\|(a_{i,i}\ve_{i,i})\|=\max_{i}|a_{i,i}|$ we may only consider matrices with
zero diagonal. Moreover, for any 
unit vectors $s,t$ we have
\begin{align*}
\Big|\sum_{i\neq j}a_{i,j}\ve_{ij}s_it_j\Big|
&\leq \|(a_{i,j})\|_\infty\sum_{i,j}\ind_{\{(i,j)\in E_A\}}\frac{1}{2}(s_i^2+t_j^2)
\\
&=
\frac{\|(a_{i,j})\|_\infty}{2}\Bigl(\sum_{i}s_i^2\sum_j\ind_{\{(i,j)\in E_A\}}
+\sum_{j}t_j^2\sum_i\ind_{\{(i,j)\in E_A\}}\Bigl)
\\
&\leq d_A\|(a_{i,j})\|_\infty.
\end{align*}
Hence, 
\begin{equation}
\label{eq:trivialbydA}
\Ex\|(a_{i,j}\ind_{\{i\neq j\}}\ve_{i,j})_{i,j}\|\leq d_A\|(a_{i,j})\|_\infty
\end{equation}
and it is enough to consider only the case $n\geq d_A\geq 3$.
\end{rem}

The proof of Theorem \ref{thm:normgraphA2} takes the most part of the paper. Here we briefly sketch the 
main ideas of this proof. Bernoulli conjecture,
formulated by  Talagrand and proven in \cite{BeL}, states that to estimate a supremum of
the Bernoulli process
one needs to decompose the index set into two parts and estimate supremum over the first part using the 
uniform bound and over the second part by the supremum of the Gaussian process. Unfortunately, there
is no algorithmic method for making such a decomposition -- a rule of thumb is that the uniform bound 
works well for large coefficients and the Gaussian bound for small ones. We try to follow this informal recipe,  
decompose vectors $s,t\in B_2^n$ into almost "flat" parts and use the uniform bound when infinity norms
of these parts are far apart. When they are of the same order we make some further technical adjustments (using 
properties of the graph $G_A$) and apply
the Gaussian bound. The crucial tool used to estimate the corresponding Gaussian process is an improvement
of van Handel's bound \cite{vH}, provided in Section 2.1.
    
We postpone the details of the proof  till the end of the paper and discuss now some consequences 
of Theorem \ref{thm:normgraphA2}.

\begin{thm}
\label{thm:uploglogdA}
For any symmetric matrix $(a_{i,j})_{i,j\leq n}$,
\[
\Ex\|(a_{i,j}\ve_{i,j})_{i,j\leq n}\|\leq 
\Log\Log(d_A)\Bigr(\max_{i}\|(a_{i,j})_j\|_2+R_A(\Log\, n)\Bigl).
\]
\end{thm}

\begin{proof}
Let $M:=\max_{i}\|(a_{i,j})_j\|_2$, $u_0=1$  and $u_k:=\exp(-(20/19)^k)$ for $k=1,2,\ldots$.
Let $k_0$ be the smallest integer such that $(\frac{20}{19})^{k_0}\geq \Log(d_A)$. Then
$k_0\sim \Log\Log(d_A)$ and $u_{k_0}\leq d_A^{-1}$. We have
\[
\Ex\|(a_{i,j}\ve_{i,j})\|\leq 
\Ex\|(a_{i,j}\ind_{\{|a_{i,j}|\leq u_{k_0}M\}}\ve_{i,j})\|
+\sum_{k=1}^{k_0}\Ex\|(a_{i,j}\ind_{\{u_kM<|a_{i,j}|\leq u_{k-1}M\}}\ve_{i,j})\|.
\]
For any $k$,
\[
d_k:=\max_{i}|\{j\colon\ |a_{i,j}|>u_kM\}|\leq u_k^{-2},
\]
so by Theorem \ref{thm:normgraphA2}
\[
\Ex\|(a_{i,j}\ind_{\{u_kM<|a_{i,j}|\leq u_{k-1}M\}}\ve_{i,j})\|
\lesssim M+R_A(\Log\, n)+d_k^{19/40}u_{k-1}M\lesssim M+R_A(\Log\, n).
\]
Moreover, using again Theorem \ref{thm:normgraphA2}
\[
\Ex\|(a_{i,j}\ind_{\{|a_{i,j}|\leq u_{k_0}M\}}\ve_{i,j})\|\lesssim 
M+R_A(\Log\, n)+d_A^{19/40}u_{k_0}M\lesssim M+R_A(\Log\, n)
\qedhere
\]

\end{proof}

Obviously, $d_A\leq n$, so Theorem \ref{thm:uploglogdA} (together with the standard 
symmetrization argument) implies that  bounds \eqref{eq:upperopnormred} and 
\eqref{eq:upperopnorm} hold up double logarithms of $n$. However, decomposing matrix into
two parts and using the Bandeira-van Handel bound one may derive conjectured upper bounds 
up to triple logarithms.

\begin{thm}
\label{thm:uplogloglogn}
For any matrix $(a_{i,j})_{i,j\leq n}$,
\[
\Ex\|(a_{i,j}\ve_{i,j})_{i,j\leq n}\|
\lesssim \Log\Log\Log\, n\Bigl(
\max_{1\leq i\leq n}\|(a_{i,j})_j\|_2
+\max_{1\leq j\leq n}\|(a_{i,j})_i\|_2+R_A(\Log\, n)\Bigr)
\]
and
\begin{align*}
\Ex\|(a_{i,j}\ve_{i,j})_{i,j\leq n}\|
\lesssim \Log\Log\Log\, n\Biggl(&
\max_{1\leq i\leq n}\|(a_{i,j})_j\|_2
+\max_{1\leq j\leq n}\|(a_{i,j})_i\|_2
\\
&
+\max_{1\leq k\leq n}\min_{I\subset [n],|I|\leq k}
\sup_{\|s\|_2,\|t\|_2\leq 1}\left\|\sum_{i,j\notin I}a_{i,j}\ve_{i,j}s_it_j\right\|_{\Log\, k}\Biggr).
\end{align*}
\end{thm}

\begin{proof}
Assume first that the matrix $(a_{i,j})$ is symmetric.
Put $M:=\max_{1\leq i\leq m}\|(a_{i,j})_j\|_2$. Estimate \eqref{eq:BvH0} yields
\[
\Ex\|(a_{i,j}\ind_{\{|a_{i,j}|\leq M\Log^{-1/2}\, n\}}\ve_{i,j})_{i,j\leq n}\|
\lesssim \max_{1\leq i\leq m}\|(a_{i,j})_j\|_2.
\]
We have
\[
\max_{i}|\{j\colon\ |a_{i,j}|>M\Log^{-1/2}\, n\}|\leq \Log\, n,
\]
hence Theorem \ref{thm:uploglogdA}, applied to a matrix 
$(a_{i,j}\ind_{\{|a_{i,j}|> M\Log^{-1/2}\, n\}})_{i,j\leq n}$ implies
\[
\Ex\|(a_{i,j}\ind_{\{|a_{i,j}|> M\Log^{-1/2}\, n\}}\ve_{i,j})_{i,j\leq n}\|
\lesssim \Log\Log\Log\, n\Bigl(\max_{1\leq i\leq m}\|(a_{i,j})_j\|_2+R_A(\Log\, n)\Bigr).
\]
Therefore, for any symmetric matrix $(a_{i,j})$,
\begin{equation}
\label{eq:symlogloglog}
\Ex\|(a_{i,j}\ve_{i,j})_{i,j\leq n}\|
\lesssim \Log\Log\Log\, n\Bigl(\max_{1\leq i\leq m}\|(a_{i,j})_j\|_2+R_A(\Log\, n)\Bigr). 
\end{equation}

Now, supppose that matrix $(a_{i,j})$ is arbitrary.
Applying \eqref{eq:symlogloglog} to the symmetric $2n\times 2n$ matrix  
$\left(\begin{array}{cc}0&A\\A^T&0\end{array}\right)$ we get the first part of the assertion.

The second part follows from the first one as in the proof of \cite[Remark 4.5]{LS}.
\end{proof}

\textbf{Organization of the paper.} In the next section we present several examples of 
applications of the main results.
In Section \ref{sec:tools} we discuss basic tools
used in the sequel, including an improvement of
the van Handel bound for norms of Gaussian matrices from \cite{vH}. In Section \ref{sec:logdA}
we derive a weaker version of Theorem \ref{thm:uploglogdA} with $\log(d_A)$ instead
of $\log\log (d_A)$ factors. The last section is devoted to the proof of Theorem \ref{thm:normgraphA2}.


\section{Examples}
\label{sec:examples}

In the first examples we discuss how to apply Corollary \ref{cor:2sided} to estimate 
$\Ex\|(\ind_E(i,j)\ve_{i,j})\|$ for various classes of graphs
$G=([n],E)$. Due to Remark \ref{rem:symmetrization} the presented  bounds 
are valid also for symmetric random matrices $\Ex\|(\ind_E(i,j)\tilde{\ve}_{i,j})_{i,j}\|$.

\begin{example}
Let $G=([n],E)$ be a graph with maximal degree $d$. Then 
\begin{equation}
\label{eq:trivial01}
\Ex\|(\ind_E(i,j)\ve_{i,j})_{i,j\leq n}\|\leq d
\end{equation}
and
\begin{equation}
\label{eq:BvH01}
\Ex\|(\ind_E(i,j)\ve_{i,j})_{i,j\leq n}\|\lesssim \sqrt{d}+\sqrt{\Log\, n}.
\end{equation}
\end{example}

The first bound follows by \eqref{eq:trivialbydA}. 
Estimate \eqref{eq:BvH01} follows from \eqref{eq:BvH0}, one may also use \eqref{eq:upper01} 
and bound the operator norm of $\ind_F$ by
the Hilbert-Schmidt norm.

\begin{example} 
Let $G=([n],E)$ be the disjoint union of $m$ complete
graphs of size $d+1$ (i.e. $n=m(d+1)$ and $\ind_E$ is a block-diagonal matrix with $m$ 
$(d+1)\times(d+1)$ blocks with ones outside diagonal).  Then
\[
\Ex\|(\ind_E(i,j)\ve_{i,j})\|\sim \sqrt{d}+\min\{d,\sqrt{\Log\, n}\}\sim
\begin{cases}
d & \mbox{if } d\leq \sqrt{\Log\, n},
\\
\sqrt{\Log\, n} & \mbox{if } \sqrt{\Log\, n}\leq d\leq \Log\, n,
\\
\sqrt{d} & \mbox{if }  d\geq \Log\, n.
\end{cases}
\]
\end{example}
Indeed, $G$ has degree $d$ and 
\[
\max_{F\subset E, |F|\leq \Log\, n}\|(\ind_{\{(i,j)\in F\}})_{i,j}\|=
\max_{F\subset [d+1]\times [d+1], |F|\leq \Log\, n}\|(\ind_{\{(i,j)\in F,i\neq j\}})_{i,j}\|
\sim \min\{d,\sqrt{\Log\, n }\}.
\]
hence the upper bound for $\Ex\|(a_{i,j}\ve_{i,j})\|$ follows by \eqref{eq:upper01}.
To derive the lower bound it is enough to consider only the case
$\Log\, n\geq d$. We then apply estimate \eqref{eq:2sided01} and
observe that
\begin{align*}
\min_{|I|\leq m-1}\max_{F\subset E, |F|\leq \Log(m-1)}
\|(\ind_{\{(i,j)\in F,i,j\notin I\}})_{i,j}\|
&=
\max_{F\subset [d]\times [d], |F|\leq \Log(m-1)}\|(\ind_{\{(i,j)\in F,i\neq j\}})_{i,j}\|
\\
&\sim \min\{d,\sqrt{\Log\, m}\}\sim \min\{d,\sqrt{\Log\, n}\}.
\end{align*}

\begin{example} 
\label{ex:largegirth}
Let $G=([n],E)$ have  maximal degree $d$ and
girth at least $\ve\Log\Log\, n$. Then $\Ex\|(\ind_E(i,j)\ve_{i,j})\|\sim_\ve \sqrt{d}$.
\end{example}

By Corollary \eqref{cor:2sided} it is enough to show that any subgraph $H=(V,F)$ of $([n],E)$
with at most $\Log\, n$ edges has spectral radius at most $C(\ve)\sqrt{d}$.
Subgraph $H$ has at most $2\Log\, n$ vertices, maximal degree at most $d$ and girth at least 
$\ve\Log\Log\, n$. Let $k=2\lfloor\frac{\ve}{2}\Log\Log\ n\rfloor$ then
$\|(\ind_F)\|\leq (\mathrm{tr}(\ind_F^k))^{1/k}$.
Since  $F$ does not contain cycles of lenghts at most $k$, we have
\[
\mathrm{tr}(\ind_F^k)=
\sum_{i_1,\ldots,i_k}\ind_{\{(i_1,i_2)\in F\}}\cdots\ind_{\{(i_{k-1},i_k)\in F\}}\ind_{\{(i_k,i_{1})\in F\}}
\leq 2\Log\, n\cdot N_{k,d}, 
\]
where $N_{k,d}\leq 2^kd^{k/2}$ is the number of closed pathes from root to itself in the $d$-regular tree.
Hence we obtain the desired bound
\[
\|(\ind_F)\|\leq (2\Log\, n)^{1/k}N_{k,d}^{1/k}\lesssim_\ve \sqrt{d}.
\]

The next example generalizes the previous one. It is close to \cite[Theorem 1.2]{MOP},
where there was a stronger assumption on the neighborhood diameter, but a more precise
bound on the operator norm.

\begin{example}
\label{ex:onecycle}
Let $G=([n],E)$  have  maximal degree $d$ and suppose that the $r$-neighborhood of every vertex
contains at most one cycle, where $r\geq \ve\Log\Log\, n$. Then 
 $\Ex\|(\ind_E(i,j)\ve_{i,j})\|\sim_\ve \sqrt{d}$.
\end{example}

Observe that all cycles of lenght at most $2r$ do not intersect. Let us remove one edge
from each cycle of length at most $2r$ and let $G_1:=([n],E_1)$, where $E_1$ contains all removed edges and 
$G_2:=([n],E_2)$, where $E_2=E\setminus E_1$. Then  $G_2$
has girth  at least $2r$ and
the maximal degree of $G_1$ is at most $1$.
Thus by \eqref{eq:trivial01} and Example \ref{ex:largegirth},
\[
\Ex\|(\ind_E(i,j)\ve_{i,j})\|\leq 
\Ex\|(\ind_{E_1}(i,j)\ve_{i,j})\|+\Ex\|(\ind_{E_2}(i,j)\ve_{i,j})\|
\lesssim_\ve 1+\sqrt{d}\lesssim \sqrt{d}.
\]

\medskip

The next example is similar to \cite[Theorem 1.16]{MOP}, where it was showned
that for a random $d$-regular graph  $G=([n],E)$, 
$\|(\ind_E(i,j)\tilde{\ve}_{i,j})\|\leq 2\sqrt{d-1}+\ve$ with probability $1-o_n(1)$.

\begin{example} 
Let $G=([n],E)$ be a random $d$-regular graph.
Then $\Ex\|(\ind_E(i,j)\ve_{i,j})\|\sim \sqrt{d}$.
\end{example}

By \eqref{eq:BvH01} we may assume that $\Log\, n\geq d$ and by \eqref{eq:trivial01}
that $d\geq 3$.
Following \cite{Bo} we say that a graph is \emph{$r$-tangle free} if the $r$-neighborhood of every vertex contains at most one cycle
By Example \ref{ex:onecycle} and \eqref{eq:trivial01} we have
\[
\Ex\|(\ind_E(i,j)\ve_{i,j})\|
\leq C\sqrt{d}\,\Pr(G \mbox{ is $\Log\Log n$-tangle free} )
+d\,\Pr(G \mbox{ is not $\Log\Log n$-tangle free}).
\]
By \cite[Lemma 27]{Bo}
\[
\Pr(\mbox{$G$ is not $\Log\Log n$-tangle free})\lesssim \frac{(d-1)^{\Log\Log n}}{n}\lesssim d^{-1/2},
\]
where the last bound follows since we assume that $\Log\, n\geq d$.

\medskip

The next example essentially recovers \cite[Corollary 1.3]{ACKM}.

\begin{example}
Let $A$ be the adjacency matrix of a $d$-regular graph and suppose that
all eigenvalues of $A$ besides the largest one are bounded in absolute value by $\lambda$.
Then $\Ex\|(a_{i,j}\ve_{i,j})\|\lesssim \lambda$. 
\end{example}

We have $\lambda \gtrsim \sqrt{d}$ by the Alon-Boppana theorem, so
by \eqref{eq:BvH01} we may assume that $\Log\, n \ge d$.
By \eqref{eq:upper01} it is enough to show that any subgraph $H=(V,F)$ of $([n],E)$
with at most $\Log\, n$ edges has spectral radius at most $C\lambda$. Subgraph $H$ has
at most $2\Log\, n$ vertices.
Let $w=n^{-1/2}(1,\ldots,1)$ be the eigenvector of $A$ corresponding to the largest eingenvalue
$d$. 
Any $v\in\er^V$ with $\|v\|_2=1$ we may represent as $v=\langle v,w\rangle w+v'$, where
$\langle v',w\rangle=0$ and $\|v'\|_2\leq 1$. Thus
\[
\|\ind_Fv\|_2\leq \|\ind_Fv'\|+\langle v,w\rangle\|\ind_Fw\|_2
\leq \lambda+\|v\|_1n^{-1/2}|V|^{1/2}dn^{-1/2}\leq \lambda+2d\frac{\Log\, n}{n}\lesssim
\lambda,
\]
where the last estimate follows since  $\Log\, n\geq d$ and $\lambda \gtrsim 1$.

\medskip

The next example shows that estimate \eqref{eq:upper01} cannot be reversed.

\begin{example}
\label{ex:minoverkneeded}
Let $1\leq d\leq n$ and $E=([d]\times [d])\cup\{(i,i)\colon\ d<i\leq n\}$
(i.e $\ind_E$ is block diagonal with one $d\times d$ block of ones and $n-d$ blocks
of single ones).
Then 
\[
\sup_{F\subset E,|F|\leq \Log\, n}\|(\ind_{\{(i,j)\in F\}})_{i,j}\|
=\sup_{F\subset [d]\times [d],|F|\leq \Log\, n}\|(\ind_{\{(i,j)\in F\}})_{i,j}\|
\sim \min\{d,\sqrt{\Log\, n}\}
\]
and the RHS of \eqref{eq:upper01} is of the order $\sqrt{d}+\min\{d,\sqrt{\Log\, n}\}$,
whereas
\[
\Ex\|(\ind_E(i,j)\ve_{i,j})_{i,j\leq n}\|=\Ex\|(\ve_{i,j})_{i,j\leq d}\|\sim \sqrt{d}.
\]
\end{example}

The last example concerns randomized circulant matrices, investigated in \cite{LS}.

\begin{example} 
Suppose that $(a_{i,j})$ is a circulant matrix, i.e. 
$a_{i,j}=b_{i-j\, \mathrm{mod}\, n}$ for a deterministic sequence $(b_i)_{i=0}^{n-1}$.
Then for any $i$ and $j$, $\|(a_{i,j})_j\|_2=\|(a_{i,j})_i\|_2=\|(b_i)_i\|_2$. Moreover,
as was shown in the proof of  \cite[Theorem 1.3]{LS}
\[
\sup_{\|s\|_2,\|t\|_2\le 1}\left\|\sum_{i,j} a_{ij}\varepsilon_{ij}s_i t_j\right\|_{\Log\, n}
\lesssim \inf_{I\subset[n],|I|\le n/4}\sup_{\|s\|_2,\|t\|_2\le 1}
\left\|\sum_{i,j\not\in I}a_{ij}\varepsilon_{ij}s_i t_j\right\|_{\Log(n/4)},
\]

Thus Theorem \ref{thm:uplogloglogn} improves \cite[Theorem 1.3]{LS} and yields that for circulant matrices
\[
\|(b_i)\|_2+R_A(\Log\, n)\lesssim \Ex\|(a_{i,j}\ve_{i,j})\|
\lesssim \Log\Log\Log\, n\big(\|(b_i)\|_2+R_A(\Log\, n)\big).
\]
\end{example}

\section{Tools}
\label{sec:tools}


We will use the following estimate for suprema of Rademachers. It is a 
special case of \cite[Lemma 5.10]{AL}.

\begin{prop}
\label{prop:supnbern}
Let $T_1,\ldots,T_n$ be nonempty bounded subsets of $\er^N$. Then
\[
\Ex\max_{k\leq n}\sup_{t\in T_k}\sum_{i=1}^N t_i\ve_i
\lesssim \max_{k\leq n}\Ex\sup_{t\in T_k}\sum_{i=1}^N t_i\ve_i
+\max_{k\leq n}\sup_{t\in T_k}\Big\|\sum_{i=1}^N t_i\ve_i\Big\|_{\Log\, n}.
\]
\end{prop}

Another useful result is the estimate on the number of connected subsets of a graph.

\begin{lem}
\label{lem:cardkconnect}
Let $H=(V_H,E_H)$ be a graph with $n_H$ vertices and maximal degree $d_H$. \\
i) For a fixed $v\in V$ the number of connected subsets  $I\subset V_H$  with cardinality $k$ containing $v$  
is at most $(4d_H)^{k-1}$.\\
ii)  The number of all connected subsets  $I\subset V_H$  with cardinality $k$  is not bigger than $n_H(4d_H)^{k-1}$. 
\end{lem}

\begin{proof}
i) We choose a connected subset $I\ni v$ by constructing its spanning tree, rooted
at $v$. In order to do it we first choose an unlabelled rooted tree with $k$ vertices and then label its vertices
by elements of $V_H$.
The number of unlabelled rooted trees is less than the number of ordered
trees with $k$ vertices, i.e., less than the $(k-1)$-th Catalan number $C_{k-1}\leq 4^{k-1}$. The root of the tree has label $v$ and the rest
of vertices may be labelled in at  most $d_H^{k-1}$ ways. \\
Part i) of the assertion immediately yields part ii).
\end{proof}






\subsection{Van Handel-type bound}

In this part we will establish the following improvement on van Handel's bound \cite{vH}.

\begin{prop}
\label{prop:notsyml2linfty}
For any $n\times m$ matrix $(a_{i,j})_{i\leq m,j\leq n}$ and $b\in (0,1]$ we have
\begin{align*}
\Ex\sup_{s\in B_2^m\cap bB_\infty^m}\sup_{t\in B_2^n\cap bB_\infty^n}
\sum_{i\leq m,j\leq n}a_{i,j}\ve_{i,j}s_it_j
&\lesssim
\max_{i}\|(a_{i,j})_j\|_2+\max_{j}\|(a_{i,j})_i\|_2
\\
&\phantom{\lesssim}+\Log((n+m)b^2)\|(a_{i,j})_{i,j}\|_\infty.
\end{align*}
\end{prop}

Let us first formulate and prove a symmetric variant  of Propostion~\ref{prop:notsyml2linfty}.

\begin{prop}
\label{prop:syml2linfty}
Let  $(\tilde{\ve}_{i,j})_{i,j}$ be a symmetric
Rademacher matrix. Then for any symmetric matrix $(a_{i,j})_{i,j\leq n}$
 and any $b\in (0,1]$,
\[
\Ex\sup_{s,t\in B_2^n\cap bB_\infty^n}\sum_{i,j\leq n}a_{i,j}\tilde{\ve}_{i,j}s_it_j
\lesssim
\max_{i}\|(a_{i,j})_j\|_2+\Log(nb^2)\|(a_{i,j})_{i,j}\|_\infty.
\]
\end{prop}

The proof uses the following, quite standard, technical lemma.

\begin{lem}
\label{lem:orderstat}
Let $Y_1,\ldots,Y_n$ be r.v's and $m_i,\sigma_i\geq 0$ be such that
\[
\Pr(|Y_i|\geq m_i+u\sigma_i)\leq e^{-u^2/2} \quad \mbox{ for every }u\geq 0 \text{ and } i=1,\ldots,n.
\]
Then
\[
\Ex\sup_{s\in B_2^n\cap bB_\infty^n}\sum_{i=1}^ns_i^2Y_i\lesssim
\max_{i}m_i+\sqrt{\Log(nb^2)}\max_{i}\sigma_i
\]
and
\[
\Ex\sup_{s\in B_2^n\cap bB_\infty^n}\sqrt{\sum_{i=1}^ns_i^2Y_i^2}\lesssim
\max_{i}m_i+\sqrt{\Log(nb^2)}\max_{i}\sigma_i.
\]
\end{lem}

\begin{proof}
Let $(x_1^*,\ldots,x_n^*)$ denote a nondecreasing rearrangement of $|x_1|,\ldots,|x_n|$. 
We set $k=n$ if $b^2\leq 1/n$, otherwise we choose
$1\leq k\leq n-1$ such that $\frac{1}{k+1}< b^2\leq \frac{1}{k}$. 
Then $\Log(nb^2)\sim \Log (n/k)$ and 
for any $s\in B_2^n\cap bB_\infty^n$,
\begin{align*}
\sum_{i=1}^ns_i^2Y_i
&\leq \sum_{i=1}^n |s_i^*|^2Y_i^*
\leq \sum_{i=1}^k |s_i^*|^2Y_i^*+\Bigl(1-\sum_{i=1}^k |s_i^*|^2\Bigr)Y_{k}^*
\\
&= \sum_{i=1}^k \Bigl(|s_i^*|^2Y_i^*+\Bigl(\frac{1}{k}-|s_i^*|^2\Bigr)Y_{k}^*\Bigr)
\leq
\frac{1}{k}(Y_1^*+\ldots+Y_k^*).
\end{align*}
For any $1\leq l\leq n$ and $u\geq 0$,
\[
\Pr(Y_l^*\geq \max_{i}m_i+u\max_i\sigma_i)\leq 
\frac{1}{l}\sum_{i=1}^n\Pr(Y_i\geq \max_{i}m_i+u\max_i\sigma_i)\leq 
\frac{n}{l}e^{-u^2/2}.
\]
Hence integration by parts yields 
$\Ex Y_l^*\leq (\Ex |Y_l^*|^2)^{1/2}
\lesssim \max_i m_i+\Log^{1/2}(n/l)\max_i\sigma_i$.
Thus
\begin{align*}
\Ex\sup_{s\in B_2^n\cap bB_\infty^n}\sum_{i=1}^ns_i^2Y_i
&\leq \frac{1}{k}\sum_{l=1}^k\Ex Y_l^*
\lesssim \max_i m_i+\frac{1}{k}\sum_{l=1}^k\sqrt{\Log\bigl(\frac{n}{l}\bigr)}\max_i\sigma_i
\\
&\lesssim \max_i m_i+\sqrt{\Log\bigl(\frac{n}{k}\bigr)}\max_i\sigma_i.
\end{align*}
Similarly,
\begin{align*}
\Ex\sup_{s\in B_2^n\cap bB_\infty^n}\sqrt{\sum_{i=1}^ns_i^2Y_i^2}
&\leq
\Ex\sqrt{\frac{1}{k}(|Y_1^*|^2+\ldots+|Y_k^*|^2)}
\leq \sqrt{\frac{1}{k}\sum_{l=1}^k\Ex |Y_l^*|^2}
\\
&\lesssim \max_i m_i+\sqrt{\frac{1}{k}\sum_{l=1}^k\Log\bigl(\frac{n}{l}\bigr)}\max_i\sigma_i
\\
&\lesssim \max_i m_i+\sqrt{\Log\bigl(\frac{n}{k}\bigr)}\max_i\sigma_i.
\qedhere
\end{align*}

\end{proof}

\begin{proof}[Proof of Proposition \ref{prop:syml2linfty}]
Let $(g_{i,j})_{i,j\leq n}$ be a symmetric Gaussian matrix (i.e., $g_{i,j}=g_{j,i}$ and
$(g_{i,j})_{i\geq j}$ are iid $\mathcal{N}(0,1)$ r.v's), independent
of $\tilde{\ve}_{i,j}$.
 We have for any matrix norm $\|\cdot\|$,
\[
\Ex\|(a_{i,j}g_{i,j})\|=\Ex\|(a_{i,j}\tilde{\ve}_{i,j}|g_{i,j}|)\|
\geq \Ex\|(a_{i,j}\tilde{\ve}_{i,j}\Ex|g_{i,j}|)\|=\sqrt{\frac{2}{\pi}}
\Ex\|(a_{i,j}\tilde{\ve}_{i,j})\|.
\]

For any symmetric matrix $B$ we have 
$\langle Bs,t\rangle
=\frac{1}{4}\bigl(\langle B(s+t),s+t\rangle-\langle B(s-t),s-t\rangle \bigr)$,
hence,
\[
\sup_{s,t\in B_2^n\cap bB_\infty^n}
\langle Bs,t\rangle
\leq 2\sup_{s\in B_2^n\cap bB_\infty^n}|\langle Bs,s\rangle|.
\]

Therefore, 
\begin{align*}
\Ex\sup_{s\in B_2^n\cap bB_\infty^n}
\sum_{i,j\leq n}a_{i,j}\tilde{\ve}_{i,j}s_it_j
&\lesssim
\Ex\sup_{s\in B_2^n\cap bB_\infty^n}\Bigl|\sum_{i,j\leq n}a_{i,j}g_{i,j}s_is_j\Bigr|
\\
&\leq \Ex\sup_{s\in B_2^n\cap bB_\infty^n}\sum_{i,j\leq n}a_{i,j}g_{i,j}s_is_j
+\Ex\sup_{s\in B_2^n\cap bB_\infty^n}\sum_{i,j\leq n}(-a_{i,j}g_{i,j}s_is_j)
\\
&=2\Ex\sup_{s\in B_2^n\cap bB_\infty^n}\sum_{i,j\leq n}a_{i,j}g_{i,j}s_is_j.
\end{align*}

Now we follow van Handel's approach from \cite{vH}. Let $g_1,g_2,\ldots,g_n$ be iid
$\mathcal{N}(0,1)$ r.v's and $Y=(Y_1,\ldots,Y_n)\sim \mathcal{N}(0,B_{-})$, where 
$B_{-}$ is the negative part of $B=(a_{i,j}^2)$. Define the new Gaussian
process $Z_{s}$ by
\[
Z_s=2\sum_{i=1}^ns_ig_i\sqrt{\sum_{j=1}^n a_{ij}^2s_j^2}+\sum_{i=1}^ns_i^2Y_i.
\]
It is shown in \cite{vH} (see the proof of Theorem 4.1 therein) 
that for any $s,s'\in \er^n$
\[
 \Ex\Bigl|\sum_{i,j\leq n}a_{i,j}g_{i,j}(s_is_j-s_i's_j')\Bigr|^2\leq \Ex|Z_s-Z_{s'}|^2.
\]
Hence the Slepian-Fernique inequality \cite[Theorem 3.15]{Ledoux-Talagrand}. yields
\[
\Ex\sup_{s\in B_2^n\cap bB_\infty^n}\sum_{i,j\leq n}a_{i,j}g_{i,j}s_is_j
\leq \Ex\sup_{s\in B_2^n\cap bB_\infty^n}Z_s.
\]
Variables $Y_i$ are centered Gaussian and (see the proof of Corollary 4.2 in \cite{vH})
$(\Ex Y_i^2)^{1/2}\leq \|(a_{i,j})_j\|_4$.
Hence Lemma \ref{lem:orderstat} applied with $m_i=0$ and $\sigma_i=\|(a_{i,j})_j\|_4$ yields
\begin{align*}
\Ex\sup_{s\in B_2^n\cap bB_\infty^n} \sum_{i=1}^ns_i^2Y_i
&\lesssim \sqrt{\Log(nb^2)}\max_{i}\|(a_{i,j})_j\|_4
\\
&\leq \sqrt{\Log(nb^2)}\max_{i,j}|a_{i,j}|^{1/2}\max_{j}\|(a_{i,j})_j\|_2^{1/2}
\\
&\leq \max_i\|(a_{i,j})_j\|_2+\Log(nb^2)\|(a_{i,j})_{i,j}\|_\infty.
\end{align*}

We have
\[
\Ex\sup_{s\in B_2^n\cap bB_\infty^n}\sum_{i=1}^ns_ig_i\sqrt{\sum_{j}a_{ij}^2s_j^2}
\leq \Ex\sup_{s\in B_2^n\cap bB_\infty^n}\sqrt{\sum_{i,j}a_{ij}^2s_j^2g_i^2}
=\Ex\sup_{s\in B_2^n\cap bB_\infty^n}\sqrt{\sum_{j}s_j^2 V_j^2},
\]
where $V_j=\sqrt{\sum_{i}a_{ij}^2g_i^2}$. The Gaussian concentration 
\cite[Lemma 3.1]{Ledoux-Talagrand} yields
\[
\Pr(|V_j|\geq \|(a_{i,j})_i\|_2+t\|(a_{i,j})_i\|_\infty)\leq e^{-t^2/2},
\]
so Lemma \ref{lem:orderstat} applied with $Y_{j}=V_j$, $m_j=\|(a_{i,j})_i\|_2$
and $\sigma_j=\|(a_{i,j})_i\|_\infty$ yields
\[
\Ex\sup_{s\in B_2^n\cap bB_\infty^n}\sum_{i=1}^ns_i\sqrt{\sum_{j}a_{ij}^2s_j^2}g_i
\lesssim  \max_{j}\|(a_{i,j})_i\|_2+\sqrt{\Log(nb^2)}\|(a_{i,j})_{i,j}\|_\infty.
\qedhere
\]

\end{proof}

\begin{proof}[Proof of Proposition \ref{prop:notsyml2linfty}]
We apply Proposition \ref{prop:syml2linfty} to the symmetric $(n+m)\times (n+m)$ matrix
$\tilde{A}=(\tilde{a}_{i,j})$ 
of the form $\tilde{A}=\left(\begin{array}{cc}0&A\\A^T&0\end{array}\right)$. 
Observe that $\max_{i,j}|\tilde{a}_{i,j}|=\max_{i,j}|a_{i,j}|$,
\[
\max_{i}\|(\tilde{a}_{i,j})_j\|_2
=\max\Bigl\{\max_{i}\|(a_{i,j})_j\|_2,\max_j\|(a_{i,j})_i\|_2\Bigr\}
\]
and
\begin{align*}
\Ex\sup_{s,t\in B_2^{n+m}\cap bB_\infty^{n+m}}
&\sum_{i,j\leq n+m}\tilde{a}_{i,j}\tilde{\ve}_{i,j}s_it_j
\geq
\Ex\sup_{s\in B_2^m\cap bB_\infty^m}\sup_{t\in B_2^n\cap bB_\infty^n}
\sum_{i\leq m,j\leq n}a_{i,j}\ve_{i,j}s_it_j.
\qedhere
\end{align*}
\end{proof}


\section{Bounds up to polylog factors} 
\label{sec:logdA}

In this section we derive weaker estimates than in Theorem \ref{thm:uploglogdA} (with powers of $\log d_A$ 
instead of $\log\log(d_A)$). They will be used in the proof of Theorem \ref{thm:normgraphA2} to estimate the
parts of Bernoulli process $(\sum_{i,j}a_{i,j}\ve_{i,j}s_it_j)_{s,t\in B_2^n}$, where coefficients $s_i$ and 
$t_j$ are of the same order. 

Let us first introduce the notation which will be used till the end of the paper. Recall that
$G_A=([n],E_A)$ is a graph associated to a given symmetric matrix 
$A=(a_{i,j})_{i,j\leq n}$.
By $\rho=\rho_A$ we denote the distance on $[n]$ induced by $E_A$. 
For $r=1,2,\ldots$ we put $G_r=G_{r}(A)=([n],E_{A,r})$, where $(i,j)\in E_{A,r}$ iff $\rho(i,j)\leq r$.
In particular $G_{1}=([n],E_A)$ and the maximal degree of $G_{r}$ is
at most $d_A+d_A(d_A-1)+\ldots+d_A(d_A-1)^{r-1}\leq d_A^r$. We say that a subset of $[n]$ is \emph{$r$-connected} if it is connected in $G_{r}$.

We denote by $\mathcal{I}(k)=\mathcal{I}(k,n)$ the family of all subsets of $[n]$ of cardinality $k$ and by $\mathcal{I}_r(k)=\mathcal{I}_r(k,A)$ the family of all $r$-connected subsets of $[n]$ of cardinality $k$.

For a set $I\subset [n]$ and a vertex $j\in [n]$ we write $I\sim_A j$ if  $(i,j)\in E_A$ for some 
$i\in I$. By $I'=I'(A)$ we denote the set of all neighbours of $I$ in $G_1$ and by $I''=I''(A)$ the set of all neighbours of $I'$ in $G_1$, i.e.,
\begin{equation}
\label{eq:defI'}
I'=\{j\in [n]\colon\ \exists_{i\in I}\ (i,j)\in E_A\},\quad
I''=\{i\in [n]\colon\ \exists_{i_0\in I,j\in [n]}\ (i_0,j), (i,j)\in E_A\}.
\end{equation}
Observe that $I$ is a subset of $I''$, but does not have to be a subset of $I'$.
Moreover $|I'|\leq d_A|I|$ and $|I''|\leq d_A^2|I|$.

By Remark \ref{rem:zerodiag} we may and will assume that $a_{i,i}=0$ for all $i$.

For $1\leq k,l\leq n$ define random variables
\[
X_{k,l}=X_{k,l}(A):=\frac{1}{\sqrt{kl}}\max_{I\in \calI(k), J\in \calI(l)}\max_{\eta_i,\eta_j'=\pm1 }
\sum_{i\in I,j\in J}a_{i,j}\ve_{i,j}\eta_i\eta_j'
\]
and their $4$-connected counteparts
\[
\overline{X}_{k,l}=\overline{X}_{k,l}(A):=\frac{1}{\sqrt{kl}}\max_{I\in \calI_4(k), J\in \calI_4(l)}
\max_{\eta_i,\eta_j'=\pm1 } \sum_{i\in I,j\in J}a_{i,j}\ve_{i,j}\eta_i\eta_j'.
\]
Set 
\begin{equation}
\label{eq:defX}
X=X(A):=\max_{1\leq k,l\leq n}X_{k,l}
=\max_{\emptyset\neq I,J\subset V}\frac{1}{\sqrt{|I||J|}}\max_{\eta_i,\eta'_j=\pm 1}
\sum_{i\in I, j\in J}a_{i,j}\ve_{i,j}\eta_{i}\eta'_j.
\end{equation}
Observe that $X$ is nonnegative.

Variables $\overline{X}_{k,l}$ are easier to estimate than $X_{k,l}$, since the number of $4$-connected subsets
is much smaller than the number of all subsets -- calculations based on this idea are made in Lemma 
\ref{lem:supoverconnA}. Lemma \ref{lem:redtosupconnA} shows that in expectation these variables do not differ
too much.


\begin{lem}
\label{lem:redtosupconnA}
For any $1\leq k,l\leq n$, 
\[
\Ex X_{k,l}\lesssim \max_{1\leq k'\leq k,1\leq l'\leq l}\Ex\overline{X}_{k',l'}+R_A(\Log(kl)).
\]
\end{lem}

\begin{proof}
Let us first fix sets $I\in \calI(k)$ and $J\in \calI(l)$.
Let $I_1,\ldots,I_r$ be  connected components of $I\cap J'$ in $G_{2}$
and $J_u:=J\cap I_u'$. Then sets $J_1,\ldots,J_r$ are disjoint.
They are also 4-connected subsets of $J$,
since otherwise there would exists a nonempty set $V\subsetneq J_u$ such that $\rho_A(V,J_u\setminus V)\geq 4$ and taking for $\tilde{V}$ neighbours of $V$ in $I_u$ we would have $\emptyset\neq \tilde{V}\subsetneq I_u$ and $\rho_A(\tilde{V},I_u\setminus \tilde{V})\geq 2$,
contradicting 2-connectivity of $I_u$.
Hence, for every  $\eta_i,\eta_j'=\pm 1$ we have
\begin{align*}
\sum_{i\in I, j\in J}a_{i,j}\ve_{i,j}\eta_i\eta_j'
&
=\sum_{i\in I\cap J', j\in J}a_{i,j}\ve_{i,j}\eta_i\eta_j'
= \sum_{u=1}^r\sum_{i\in I_u,j\in J_u}a_{i,j}\ve_{i,j}\eta_i\eta_j'
\leq \sum_{u=1}^r\overline{X}_{|I_u|,|J_u|}\sqrt{|I_u||J_u|}
\\
&\leq
\max_{k'\leq k,l'\leq l}\overline{X}_{k',l'}\sum_{u=1}^r\sqrt{|I_u||J_u|}
\leq 
\max_{k'\leq k,l'\leq l}\overline{X}_{k',l'}\Bigl(\sum_{u=1}^r|I_u|\Bigr)^{1/2}\Bigl(\sum_{u=1}^r|J_u|\Bigr)^{1/2}
\\
&\leq 
\max_{k'\leq k,l'\leq l}\overline{X}_{k',l'}\sqrt{|I||J|}.
\end{align*}
 
Taking the supremum over all sets  $I\in \calI(k)$, $J\in \calI(l)$ and
 $\eta_i,\eta_j'=\pm 1$ we get
\begin{equation}
\label{eq:flatvsflatconnA}
X_{k,l}
\leq \max_{k'\leq k,l'\leq l}\overline{X}_{k',l'}.
\end{equation}
Observe that 
\begin{align*}
\max_{I\in \calI_4(k'),J\in \calI_4(l')}\max_{\eta_i,\eta'_j=\pm 1 }\frac{1}{\sqrt{k'l'}}
\Bigl\|\sum_{i\in I,j\in J}a_{i,j}\ve_{i,j}\eta_i\eta'_j\Bigr\|_{\Log(kl)}
&\leq \sup_{\|s\|_2,\|t\|_2\leq 1}\Bigl\|\sum_{i,j}a_{i,j}\ve_{i,j}s_it_j\Bigr\|_{\Log(kl)}
\\
&=R_A(\Log(kl)).
\end{align*}
Thus, by Proposition \ref{prop:supnbern}
\[
\Ex\max_{k'\leq k,l'\leq l}\overline{X}_{k',l'}
\lesssim  \max_{k'\leq k,l'\leq l}\Ex \overline{X}_{k',l'}+R_A(\Log(kl)).
\]
\end{proof}

\begin{lem}
\label{lem:supoverconnA}
We have for any $1\leq k,l\leq n$,
\[
\Ex \overline{X}_{k,l}\lesssim \sqrt{\Log\, d_A}\max_{i}\|(a_{i,j})_j\|_2+R_A(\Log\, n).
\]
\end{lem}

\begin{proof}
Obviously $\overline{X}_{k',l'}\leq \|(a_{i,j}\ve_{i,j})\|$, so by \eqref{eq:trivialbydA}
we may assume that $n\geq d_A\geq 3$. By the symmetry it is enough to consider only the case $l\geq k$.
By Lemma \ref{lem:cardkconnect},
$2^k|\calI_4(k)|\leq n(8d_A^4)^k\leq nd_A^{6k}$.

We have
\[
\overline{X}_{k,l}\leq
\frac{1}{\sqrt{k}}\max_{I\in \mathcal{I}_4(k)}\max_{\eta_i=\pm 1}
\sup_{\|t\|_2\leq 1}\sum_{i\in I}\sum_{j}a_{i,j}\ve_{i,j}\eta_it_j.
\]
For any fixed $I\in \mathcal{I}_4(k)$ and $\eta_i=\pm 1$,
\begin{align*}
\Ex\frac{1}{\sqrt{k}}\sup_{\|t\|_2\leq 1}\sum_{i\in I}\sum_{j}a_{i,j}\ve_{i,j}\eta_it_j
&=\frac{1}{\sqrt{k}}\Ex\Bigl(\sum_{j}\Bigr(\sum_{i\in I}a_{i,j}\ve_{i,j}\eta_i\Bigl)^2\Bigr)^{1/2}
\\
&\leq \frac{1}{\sqrt{k}}\Bigl(\sum_{j}\Ex\Bigr(\sum_{i\in I}a_{i,j}\ve_{i,j}\eta_i\Bigl)^2\Bigr)^{1/2}
=\frac{1}{\sqrt{k}}\Bigl(\sum_{j}\sum_{i\in I}a_{i,j}^2\Bigr)^{1/2}
\\
&=\frac{1}{\sqrt{k}}\Bigl(\sum_{i\in I}\sum_{j}a_{i,j}^2\Bigr)^{1/2}
\leq \max_{i}\|(a_{i,j})_j\|_2.
\end{align*}

In the case $n\geq d_A^{6k}$, $\log(2^k|\calI_4(k)|)\lesssim \log n$ and
\[
\max_{I\in \mathcal{I}_4(k)}\max_{\eta_i=\pm 1}\sup_{\|t\|_2\leq 1}\frac{1}{\sqrt{k}}
\Bigr\|\sum_{i\in I}\sum_{j}a_{i,j}\ve_{i,j}\eta_it_j \Bigl\|_{\Log(2^k|\calI_4(k)|)}
\lesssim R_A(\log n).
\]
In the case $n\leq d_A^{6k}$ we have $\log(2^k|\calI_4(k)|)\lesssim k\log d_A$ and
\begin{align*}
\max_{I\in \mathcal{I}_4(k)}&\max_{\eta_i=\pm 1}\sup_{\|t\|_2\leq 1}\frac{1}{\sqrt{k}}
\Bigr\|\sum_{i\in I}\sum_{j}a_{i,j}\ve_{i,j}\eta_it_j \Bigl\|_{\Log(2^k|\calI_4(k)|)}
\\
&\lesssim \max_{I\in \mathcal{I}_4(k)}\max_{\eta_i=\pm 1}\sup_{\|t\|_2\leq 1}\sqrt{\log d_A}
\Bigl(\sum_{i\in I}\sum_{j}(a_{i,j}\eta_it_j)^2\Bigr)^{1/2}
\\
&= \max_{I\in \mathcal{I}_4(k)}\max_j \sqrt{\log d_A}\Bigl(\sum_{i\in I}a_{i,j}^2\Bigr)^{1/2}
\leq \sqrt{\log d_A}\max_{j}\|(a_{i,j})_i\|_2.
\end{align*}
The assertion follows by Proposition \ref{prop:supnbern}.
\end{proof}

\begin{cor}
\label{cor:estsupflatA}
We have
\[
\Ex X=
\Ex\max_{\emptyset\neq I,J\subset V}\frac{1}{\sqrt{|I||J|}}\max_{\eta_i,\eta_j'=\pm 1}
\sum_{i\in I, j\in J}a_{i,j}\ve_{i,j}\eta_{i}\eta'_j
\lesssim \sqrt{\Log\, d_A}\max_{i}\|(a_{i,j})_j\|_2+R_A(\Log\, n).
\]
\end{cor}

\begin{proof}
Lemmas \ref{lem:redtosupconnA} and \ref{lem:supoverconnA} imply that for a fixed 
$1\leq k,l\leq n$
\[
\Ex X_{k,l}\lesssim \max_{k'\leq k,l'\leq l} \Ex\overline{X}_{k',l'}+R_A(\Log(kl)) 
\lesssim \sqrt{\Log\, d_A}\max_{i}\|(a_{i,j})_j\|_2+R_A(\Log\, n).
\]
Moreover,
\begin{align*}
\max_{1\leq k,l\leq n}\max_{I\in \calI(k),J\in \calI(l)}\frac{1}{\sqrt{kl}}\max_{\eta_i,\eta_j'=\pm 1}
&\Big\|\sum_{i\in I, j\in J}a_{i,j}\ve_{i,j}\eta_{i}\eta_j\Big\|_{\Log(n^2)}
\\
&\leq
\sup_{\|t\|_2\leq 1,\|s\|_2\leq 1}\Big\|\sum_{i,j\in V}a_{i,j}\ve_{i,j}t_{i}s_j\Big\|_{2\Log\, n}
\lesssim R_A(\Log\, n)
\end{align*}
and the assertion follows by Proposition \ref{prop:supnbern}.
\end{proof}

\begin{prop}
\label{prop:normgraphA1}
For any symmeric matrix $(a_{i,j})_{i,j\leq n}$  we have
\[
\Ex\bigl\|(a_{i,j}\ve_{i,j})_{i,j\leq n}\bigr\|
\lesssim \Log^{3/2}(d_A) \max_{i}\|(a_{i,j})_j\|_2+ \Log(d_A) R_A(\Log\, n).
\]
\end{prop}

\begin{proof}
By Remark \ref{rem:zerodiag} we may assume that $a_{i,i}=0$ for all $i$ and  
$n\geq d_A\geq 3$.

For vectors $s,t$ and integers $k,l$ we define sets
\[
I_k(s)=\{i\leq n\colon\ e^{-k-1}<|s_i|\leq e^{-k}\},\quad
J_l(t)=\{j\leq n\colon\ e^{-l-1}<|t_j|\leq e^{-l}\}.
\] 
Observe that for any $s,t,k,l$,
\[
\sum_{i\in I_k(s),j\in J_l(t)}a_{i,j}\ve_{i,j}s_it_j
\leq e^{-k-l}\max_{\eta_i,\eta'_j=\pm 1}\sum_{i\in I_k(s), j\in J_l(t)}a_{i,j}\ve_{i,j}\eta_{i}\eta'_{j},
\]
therefore
\[
\bigl\|(a_{i,j}\ve_{i,j})_{i,j\leq n}\bigr\|\leq 
\sup_{\|s\|_2,\|t\|_2\leq 1}\sum_{k,l} e^{-k-l}\sup_{\eta_i,\eta_j'=\pm 1}
\sum_{i\in I_k(s), j\in J_l(t)}a_{i,j}\ve_{i,j}\eta_{i}\eta'_{j}.
\]
We have
\begin{align*}
\sum_{k}\sum_{l\geq k+\log d_A} &e^{-k-l}\max_{\eta_i,\eta_j'=\pm 1}
\sum_{i\in I_k(s), j\in J_l(t)}a_{i,j}\ve_{i,j}\eta_{i}\eta'_{j}
\\
&\leq \sum_{k}\sum_{l\geq k+\log d_A} e^{-k} \sum_{i\in I_k(s),j\in J_l(t)}|a_{i,j}|e^{-l}
\\
&= \sum_{k}e^{-k}\sum_{i\in I_k(s)}\sum_{j}|a_{i,j}|\sum_{l\geq k+\log d_A}
e^{-l}\ind_{\{j\in J_l(t)\}}
\\
&\leq \sum_{k}\sum_{i\in I_k(s)}\sum_{j}|a_{i,j}|e^{-2k-\log d_A}
\leq \|(a_{i,j})\|_\infty\sum_{k}\sum_{i\in I_k(s)}e^{-2k}
\\
&\leq \|(a_{i,j})\|_\infty \sum_{k}\sum_{i\in I_k(s)}e^2s_i^2=e^2\|s\|_2^2\|(a_{i,j})\|_\infty.
\end{align*}
In the same way we show that
\[
\sum_{l}\sum_{k\geq l+\log d_A} e^{-k-l}\max_{\eta_i,\eta_j'=\pm 1}
\sum_{i\in I_k(s), j\in J_l(t)}a_{i,j}\ve_{i,j}\eta_{i}\eta'_{j}
\leq e^2\|t\|_2^2\|(a_{i,j})\|_\infty.
\]
Moreover, for any $s,t$,
\begin{align*}
\sum_{k,l\colon |k-l|<\log d_A} e^{-k-l}\sup_{\eta_i,\eta_j'=\pm 1}
\sum_{i\in I_k(s), j\in J_l(t)}
&a_{i,j}\ve_{i,j}\eta_{i}\eta'_{j}
\\
&\leq X\sum_{k,l\colon |k-l|<\log d_A}^\infty e^{-k-l}\sqrt{|I_k(s)||J_l(t)|}.
\end{align*}
For any fixed integer $r$
\begin{align*}
\sum_{k} e^{-k-(k+r)}\sqrt{|I_k(s)||J_{k+r}(t)|}
&\leq 
\Bigl(\sum_{k} e^{-2k}|I_k(s)|\Bigr)^{1/2}\Bigl(\sum_{k} e^{-2(k+r)}|J_{k+r}(t)|\Bigr)^{1/2}
\\
&\leq e^2\|t\|_2\|s\|_2.
\end{align*}
Hence,
\begin{align*}
\bigl\|(a_{i,j}\ve_{i,j})_{i,j\leq n}\bigr\|
&\leq \sup_{\|s\|_2,\|t\|_2\leq 1}e^2((\|s\|^2+\|t\|^2)\|(a_{i,j})\|_\infty+(2\log d_A+1) X\|t\|_2\|s\|_2)
\\
&\leq e^2(2\|(a_{i,j})\|_\infty+(2\log d_A+1)X)
\end{align*}
and the assertion follows by Corollary \ref{cor:estsupflatA}.
\end{proof}


\section{Proof of Theorem \ref{thm:normgraphA2}}

By 
Remark \ref{rem:zerodiag} we may assume that $a_{i,i}=0$ for all $i$ and $n\geq d_A\geq 3$.

For $k=1,2,\ldots$ and $t,s\in B_2^n$ we define
\[
I_k(s):=\{i\colon\  d_A^{-k/40}<|s_i|\leq d_A^{(1-k)/40}\},\quad
J_l(t):=\{j\colon\  d_A^{-l/40}<|t_j|\leq d_A^{(1-l)/40}\}.
\]
Then
\begin{equation}
\label{eq:sumJlA}
\sum_{k\geq 1}d_A^{-k/20}|I_k(s)|\leq \|s\|_2^2,\quad
\sum_{l\geq 1}d_A^{-l/20}|J_l(t)|\leq \|t\|_2^2
\end{equation}
and
\[
\Ex\|(a_{i,j}\ve_{i,j})_{i,j\leq n}\|=
\Ex\sup_{\|s\|_2\leq 1}\sup_{\|t\|_2\leq 1}
\sum_{k,l\geq 1}\sum_{i\in I_k(s)}\sum_{j\in J_l(t)}a_{i,j}\ve_{i,j}s_it_j.
\]

Observe that for any $s,t\in B_2^n$,
\begin{align*}
\Bigl|\sum_{k\geq 1}\sum_{l\geq k+41}\sum_{i\in I_k(s)}\sum_{j\in J_l(t)}
a_{i,j}\ve_{i,j}s_it_j\Bigr|
&\leq \sum_{k\geq 1}\sum_{i\in I_k(s)}|s_i|\sum_{l\geq k+41}\sum_{j\in J_l(t)}
|a_{i,j}||t_j|
\\
&\leq \sum_{k\geq 1}\sum_{i\in I_k(s)}|s_i|d_A^{-(k+40)/40}\sum_{j}|a_{i,j}|
\\
&\leq \|(a_{i,j})\|_\infty \sum_{k\geq 1}\sum_{i\in I_k(s)}s_i^2\leq \|(a_{i,j})\|_\infty.
\end{align*}
Similarily,
\[
\Bigl|\sum_{l\geq 1}\sum_{k\geq l+41}\sum_{i\in I_k(s)}\sum_{j\in J_l(t)}
a_{i,j}\ve_{i,j}s_it_j\Bigr|\leq \|(a_{i,j})\|_\infty.
\]
Hence it is enough  to estimate
\[
\sum_{r=-40}^{40}\Ex \sup_{\|s\|_2\leq}\sup_{\|t\|_2\leq 1}
\sum_{\substack{k,l\geq 1\\l-k=r}}\sum_{i\in I_k(s)}\sum_{j\in J_l(t)}a_{i,j}\ve_{i,j}s_it_j.
\]
By symmetry it is enough to bound only the terms with $r\geq 0$. Let $X$ be defined by
\eqref{eq:defX}.

Then for a fixed $r\geq 0$ and $\alpha>0$,
\begin{align*}
\sup_{\|s\|_2\leq 1}&\sup_{\|t\|_2\leq 1}
\sum_{k\geq 1}\sum_{i\in I_k(s)}\sum_{j\in J_{k+r}(t)}
\ind_{\{|J_{k+r}(t)|\geq \alpha|I_k(s)|\}} a_{i,j}\ve_{i,j}s_it_j
\\
&\leq
\sup_{\|s\|_2\leq 1}\sup_{\|t\|_2\leq 1}
\max_{\eta_i,\eta_j'=\pm 1}\sum_{k\geq 1}\sum_{i\in I_k(s)}\sum_{j\in J_{k+r}(t)}
d_A^{(2-2k-r)/40}\ind_{\{|J_{k+r}(t)|\geq \alpha|I_k(s)|\}} a_{i,j}\ve_{i,j}\eta_i\eta_j'
\\
&\leq X\sup_{\|s\|_2\leq 1}\sup_{\|t\|_2\leq 1}\sum_{k\geq 1}
d_A^{(2-2k-r)/40}\sqrt{|I_k(s)||J_{k+r}(t)|}\ind_{\{|J_{k+r}(t)|\geq \alpha|I_k(s)|\}}
\\
&\leq \alpha^{-1/2}X\sup_{\|s\|_2\leq 1}\sup_{\|t\|_2\leq 1}\sum_{k\geq 1} 
d_A^{(2-2k-r)/40}|J_{k+r}(t)|
\leq \alpha^{-1/2}d_A^{(r+2)/40}X,
\end{align*}
where the last inequality follows by \eqref{eq:sumJlA}.

Hence Corollary \ref{cor:estsupflatA} yields
\begin{align*}
\Ex\sup_{\|s\|_2\leq 1}\sup_{\|t\|_2\leq 1}
\sum_{k\geq 1}\sum_{i\in I_k(s)}\sum_{j\in J_{k+r}(t)}
\ind_{\{|J_{k+r}(t)|\geq d_A^{(2r+5)/40}|I_k(s)|\}}& a_{i,j}\ve_{i,j}s_it_j
\\
&\lesssim \max_i\|(a_{i,j})_j\|_2+R_A(\log n).
\end{align*}
In a similar way we show that
\begin{align*}
\sup_{\|s\|_2\leq 1}\sup_{\|t\|_2\leq 1}\sum_{k\geq 1}\sum_{i\in I_k(s)}
&\sum_{j\in J_{k+r}(t)} 
\ind_{\{|J_{k+r}(t)|\leq \alpha|I_k(s)|\}} a_{i,j}\ve_{i,j}s_it_j
\\
&\leq 
\alpha^{1/2}X\sup_{\|s\|_2\leq 1}\sup_{\|t\|_2\leq 1}\sum_{k\geq 1} 
d_A^{(2-2k-r)/40}|I_{k}(s)|
\leq \alpha^{1/2}d_A^{(2-r)/40}X
\end{align*}
and
\begin{align*}
\Ex\sup_{\|s\|_2\leq 1}\sup_{\|t\|_2\leq 1}
\sum_{k\geq 1}\sum_{i\in I_k(s)}\sum_{j\in J_{k+r}(t)}
\ind_{\{|J_{k+r}(t)|\leq d_A^{(2r-5)/40}|I_k(s)|\}}
&a_{i,j}\ve_{i,j}s_it_j
\\
&\lesssim \max_i\|(a_{i,j})_j\|_2+R_A(\log n).
\end{align*}
Hence it is enough to bound for $r=0,1,\ldots,40$,
\[
\Ex\sup_{\|s\|_2\leq 1}\sup_{\|t\|_2\leq 1}
\sum_{k\geq 1}\sum_{i\in I_k(s)}\sum_{j\in J_{k+r}(t)}
\ind_{\left\{d_A^{(2r-5)/40}<|J_{k+r}(t)|/|I_k(s)|< d_A^{(2r+5)/40}\right\}} a_{i,j}\ve_{i,j}s_it_j.
\]

Recall definition \eqref{eq:defI'} of sets $I'$ and $I''$.
Let us fix $0\leq r\leq 40$ and $k,t,s$ such that
$d_A^{(2r-5)/40}<|J_{k+r}(t)|/|I_k(s)|< d_A^{(2r+5)/40}$. Let $|I_k(s)|=m$.
Let us consider the following greedy algorithm with output being a subset $\{i_1,\ldots,i_M\}$ of $I_k(s)$ of size $M\leq m$
\begin{itemize}
\item	In the first step we pick a vertex $i_1\in I_k(s)$ with maximal number of neighbours in $J_{k+r}(t)$.
\item  Once we have $\{i_1,\ldots,i_N\}$ and $N<M$, we pick  $i_{N+1}\in I_k(s)\setminus \{i_1,\ldots,i_N\}$ with maximal number of neighbours in  $J_{k+r}(t)\setminus \{i_1,\ldots,i_N\}'$.
\end{itemize}
If $l_N$ is the number of neighbours of $i_N$ in $J_{k+r}(t)\setminus \{i_1,\ldots,i_{N-1}\}'$, then $l_1\ge l_2\ge \ldots \ge l_M$, so $Ml_M \le |J_{k+r}(t)|$.
Choose $M\leq m$ to be the maximal integer so that $l_M\geq d_A^{(r+18)/40}$ and set
$I:=\{i_1,\ldots,i_M\}$. This way we construct a subset $I\subset I_k(s)$  with
cardinality $|I|\leq d_A^{-(r+18)/40}|J_{k+r}(t)|\leq d_A^{(r-13)/40}m$ such that
for every $i\in I_k(s)\setminus I$, 
$|\{ j\in J_{k+r}(t)\setminus I'\colon i\sim_A j \}|<d_A^{(r+18)/40}$. 
  Note that if $(i,j)\in E_A$ and $(i,j)\in (I_k(s)\times J_{k+r}(t))\setminus (I''\times I')$, then $j\in J_{k+r}(t)\setminus I'$. Therefore,
\begin{align}
\notag
\sum_{(i,j)\in  I_k(s)\times J_{k+r}(t))\setminus (I''\times I')}  |a_{i,j}||s_it_j|
&\leq \|(a_{i,j})\|_\infty\sum_{i\in I_k(s)}|s_i|d_A^{(r+18)/40}d_A^{(1-k-r)/40}
\\
\label{eq:estIkJr1A}
&\leq \|(a_{i,j})\|_\infty d_A^{19/40}\sum_{i\in I_k(s)}s_i^2.
\end{align}

Let 
\[
s'=(s'_i)_{i\in I''\cap I_k(s)},\quad t'=(t'_j)_{j\in I'\cap J_{k+r}(s)},
\]
where
\[
s'_i:=\frac{s_i}{\|(s_i)_{i\in I_k(s)}\|_2},\quad 
t'_j:=\frac{t_j}{\|(t_j)_{j\in J_{k+r}(t)}\|_2}.
\]
Then 
\begin{align*}
&\|s'\|_2\leq 1,\quad \|s'\|_\infty\leq d_A^{1/40}|I_k(s)|^{-1/2}=d_A^{1/40}m^{-1/2},
\\
&\|t'\|_2\leq 1,\quad 
\|t'\|_\infty\leq d_A^{1/40}|J_{k+r}(t)|^{-1/2}\leq d_A^{(7-2r)/80}m^{-1/2}.
\end{align*}
Hence, 
\begin{equation}
\label{eq:estIkJr2A}
\sum_{(i,j)\in (I_k(s)\times J_{k+r}(t))\cap (I''\times I')}\ve_{i,j}a_{i,j}s_it_j
\leq Y_{m,r}\|(s_i)_{i\in I_k(s)}\|_2\|(t_j)_{j\in J_{k+r}(t)}\|_2,
\end{equation} 
where
\[
Y_{m,r}:=\max_{|I|\leq d_A^{(r-13)/40}m}
\sup_{s\in B_2^n\cap d_A^{1/40}m^{-1/2}B_\infty^n}\
\sup_{t\in B_2^n\cap d_A^{(7-2r)/80}m^{-1/2}B_\infty^n}\
\sum_{i\in I'',j\in I'}a_{i,j}\ve_{i,j}s_it_j.
\]

Define $Y_r:=\max_{1\leq m\leq n}Y_{m,r}$. 
The main advantage of introducing the variables $Y_r$ is that in their definition, 
the suprema are taken over delocalized vectors from the sphere of dimension $|I|$ with all
coordinates bounded essentially by $d_A^C|I|^{-1/2}$, and Proposition~\ref{prop:notsyml2linfty} 
can be used to estimate these quantities.

Estimates \eqref{eq:estIkJr1A} 
and \eqref{eq:estIkJr2A} yield
\begin{align*}
\Ex&\sup_{\|s\|_2\leq 1}\sup_{\|t\|_2\leq 1}
\sum_{k\geq 1}\sum_{i\in I_k(s)}\sum_{j\in J_{k+r}(t)}
\ind_{\Big\{d_A^{(2r-5)/40}<|J_{k+r}(t)|/|I_k(s)|< d_A^{(2r+5)/40}\Big\}} a_{i,j}\ve_{i,j}s_it_j
\\
&\leq \sup_{\|s\|_2\leq 1}\sup_{\|t\|_2\leq 1}\Biggl(
\sum_{k\geq 1}d_A^{19/40}\|(a_{i,j})\|_\infty\sum_{i\in I_k(s)}s_i^2
+\Ex Y_r\sum_{k\geq 1}
\|(s_i)_{i\in I_k(s)}\|_2\|(t_j)_{j\in J_{k+r}(t)}\|_2\Biggr)
\\
&\leq d_A^{19/40}\|(a_{i,j})\|_\infty
+\Ex Y_r\sup_{\|s\|_2\leq 1}\Bigl(\sum_{k\geq 1}\|(s_i)_{i\in I_k(s)}\|_2^2\Bigr)^{1/2}
\sup_{\|t\|_2\leq 1}\Bigl(\sum_{k\geq 1}\|(t_j)_{j\in J_{k+r}(t)}\|_2^2\Bigr)^{1/2}
\\
&\leq d_A^{19/40}\|(a_{i,j})\|_\infty+\Ex Y_r.
\end{align*}

Therefore, to establish Theorem \ref{thm:normgraphA2} it is enough to prove the following lemma.

\begin{lem}
\label{lem:maxYmrA}
For every $0\leq r\leq 40$,
\[
\Ex Y_r=\Ex\max_{1\leq m\leq n}Y_{m,r}
\lesssim \max_{i}\|(a_{i,j})_j\|_2+R_A(\log n)+\Log(d_A)\|(a_{i,j})\|_\infty.
\]
\end{lem}

First we show a connected counterpart to Lemma \ref{lem:maxYmrA}.

\begin{lem}
\label{lem:maxconnYA}
We have
\begin{align*}
\Ex \max_{1\leq k\leq n}\max_{I\in \mathcal{I}_4(k)}\
\sup_{s\in B_2^n\cap d_A^{3/8}k^{-1/2}B_\infty^n}\
&\sup_{t\in B_2^n\cap d_A^{-1/16}k^{-1/2}B_\infty^n}\
\sum_{i\in I'',j\in I'}a_{i,j}\ve_{i,j}s_it_j
\\
&\lesssim \max_{i}\|(a_{i,j})_j\|_2+R_A(\log n)+\Log(d_A)\|(a_{i,j})\|_\infty.
\end{align*}
\end{lem}

\begin{proof}
Let us first fix $k$ and $I\in \mathcal{I}_4(k)$. Then $|I'|\leq d_Ak$ and $|I''|\leq d_A^2k$. 
Proposition~\ref{prop:notsyml2linfty}, applied with 
$(a_{i,j})=(a_{i,j})_{i\in I'',j\in I'}$, $n=|I''|$, $m=|I'|$, and $b=d_A^{3/8}k^{-1/2}$ yields
\begin{align*}
\Ex \sup_{s\in B_2^n\cap d_A^{3/8}k^{-1/2}B_\infty^n}\
\sup_{t\in B_2^n\cap d_A^{-1/16}k^{-1/2}B_\infty^n}\
\sum_{i\in I'',j\in I'}&a_{i,j}\ve_{i,j}s_it_j
\\
&\lesssim \max_{i}\|(a_{i,j})_j\|_2+\Log(d_A)\|(a_{i,j})\|_\infty.
\end{align*}

By Lemma \ref{lem:cardkconnect}, $|\mathcal{I}_4(k)|\leq n(4d_A^4)^k\leq \max\{n^2,d_A^{12k}\}$
(recall that we assume that $d_A\geq 3$).
We have
\begin{align*}
&\sup_{s\in B_2^n\cap d_A^{3/8}k^{-1/2}B_\infty^n}\
\sup_{t\in B_2^n\cap d_A^{-1/16}k^{-1/2}B_\infty^n}
\Bigl\|\sum_{i,j}a_{i,j}\ve_{i,j}s_it_j\Bigr\|_{\Log(|\mathcal{I}_4(k)|)}
\\
&\leq 
\sup_{s,t\in B_2^n}
\Bigl\|\sum_{i,j}a_{i,j}\ve_{i,j}s_it_j\Bigr\|_{2\Log(n)}
+\sup_{s\in B_2^n}\sup_{t\in B_2^n\cap d_A^{-1/16}k^{-1/2}B_\infty^n}
\Bigl\|\sum_{i,j}a_{i,j}\ve_{i,j}s_it_j\Bigr\|_{12k\Log(d_A)}
\\
&\leq R_A(2\Log n)+\sup_{s\in B_2^n}\sup_{t\in B_2^n\cap d_A^{-1/16}k^{-1/2}B_\infty^n}
\sqrt{12k\Log(d_A)}
\Bigl(\sum_{i,j}a_{i,j}^2s_i^2t_j^2\Bigr)^{1/2}
\\
&\lesssim R_A(\Log\, n)+
\sqrt{k\Log(d_A)}\max_{i}d_A^{-1/16}k^{-1/2}\Bigl(\sum_{j}a_{i,j}^2\Bigr)^{1/2}
\\
&\leq R_A(\Log\, n)+ \sqrt{\Log(d_A)}d_A^{-1/16}\max_{i}\|(a_{i,j})_j\|_2\lesssim 
R_A(\log n)+ \max_{i}\|(a_{i,j})_j\|_2.
\end{align*}
Hence, by Proposition \ref{prop:supnbern},
\begin{align*}
\Ex \max_{I\in \mathcal{I}_4(k)}
\sup_{s\in B_2^n\cap d_A^{3/8}k^{-1/2}B_\infty^n}\
&\sup_{t\in B_2^n\cap d_A^{-1/16}k^{-1/2}B_\infty^n}\
\sum_{i\in I'',j\in I'}a_{i,j}\ve_{i,j}s_it_j
\\
&\lesssim \max_{i}\|(a_{i,j})_j\|_2+R_A(\log n)+\Log(d_A)\|(a_{i,j})\|_\infty.
\end{align*}
Applying again Proposition \ref{prop:supnbern} and observing that
\[
\max_{1\leq k\leq n}
\sup_{s\in B_2^n\cap d_A^{3/8}k^{-1/2}B_\infty^n}\
\sup_{t\in B_2^n\cap d_A^{-1/16}k^{-1/2}B_\infty^n}\
\Bigl\|\sum_{i,j}a_{i,j}\ve_{i,j}s_it_j\Bigr\|_{\Log(n)}
\leq R_A(\log n)
\]
we get the assertion.
\end{proof}

\begin{proof}[Proof of Lemma \ref{lem:maxYmrA}]
Let
\[
Z_k:=\max_{I\in \mathcal{I}_4(k)}
\sup_{s\in B_2^n\cap d_A^{3/8}k^{-1/2}B_\infty^n}\
\sup_{t\in B_2^n\cap d_A^{-1/16}k^{-1/2}B_\infty^n}\
\sum_{i\in I'',j\in I'}a_{i,j}\ve_{i,j}s_it_j.
\]

Let us fix $0\leq r\leq 40$, $I\subset V$ such that $|I|\leq d_A^{(r-13)/40}m$, 
$s\in B_2^n\cap d_A^{1/40}m^{-1/2}B_\infty^n$ and $t\in B_2^n\cap d_A^{(7-2r)/80}m^{-1/2}B_\infty^n$.
Let $I_1,\ldots,I_l$ be 4-connected
components of $I$. Then $\{I_1',\ldots,I_l'\}$ is a partition of $I'$,
 $\{I_1'',\ldots,I_l''\}$ is a partition of $I''$ and
\begin{equation}
\label{eq:estYkmA0}
\sum_{i\in I'',j\in I'}a_{i,j}\ve_{i,j}s_it_j
=\sum_{u=1}^l\sum_{i\in I_u'',j\in I_u'}a_{i,j}\ve_{i,j}s_it_j.
\end{equation}
Define 
\begin{align*}
U_1&:=\bigl\{1\leq u\leq l\colon\
\|(s_i)_{i\in I_{u}''}\|_2\geq d_A^{(12-r)/80}m^{-1/2}\sqrt{|I_u|}\bigr\},
\\
U_2&:=\bigl\{1\leq u\leq l\colon\
 \|(t_j)_{j\in I_{u}'}\|_2\geq d_A^{(12-r)/80}m^{-1/2}\sqrt{|I_u|}\bigr\}.
\end{align*}
For  $u\in U_1\cap U_2$ define vectors
\[
\tilde{s}(u):=\frac{(s_i)_{i\in I_u''}}{\|(s_i)_{i\in I_{u}''}\|_2},
\quad \tilde{t}(u):=\frac{(t_j)_{j\in I_u'}}{\|(t_j)_{j\in I_{u}'}\|_2}.
\]
Then $\|\tilde{s}(u)\|_2=\|\tilde{t}(u)\|_2=1$, 
\begin{align*}
\|\tilde{s}(u)\|_\infty&\leq d_A^{(r-12)/80}m^{1/2}|I_u|^{-1/2}\|s\|_\infty
\leq d_A^{(r-10)/80}|I_u|^{-1/2}\leq d_A^{3/8}|I_u|^{-1/2},
\\
\|\tilde{t}(u)\|_\infty&\leq d_A^{(r-12)/80}m^{1/2}|I_u|^{-1/2}\|t\|_\infty
\leq d_A^{-(r+5)/80}|I_u|^{-1/2}\leq d_A^{-1/16}|I_u|^{-1/2}.
\end{align*}
Hence
\begin{align}
\notag
\sum_{u\in U_1\cap U_2}\sum_{i\in I_u'',j\in I_u'}a_{i,j}\ve_{i,j}s_it_j
&\leq \sum_{u\in U_1\cap U_2}Z_{|I_u|}\|(s_i)_{i\in I_{u}''}\|_2\|(t_j)_{j\in I_{u}'}\|_2
\\
\notag
&\leq\max_{k}Z_k\Bigl(\sum_{u\leq l}\|(s_i)_{i\in I_{u}''}\|_2^2\Bigr)^{1/2}
\Bigl(\sum_{u\leq l}\|(t_j)_{j\in I_{u}'}\|_2^2\Bigr)^{1/2}
\\
\label{eq:estYkmA1}
&\leq \max_{k}Z_k.
\end{align}

Observe that
\[
\sum_{u\notin  U_1}\|(s_i)_{i\in I_{u}''}\|_2^2
\leq \sum_{u} d_A^{(12-r)/40}m^{-1}|I_u|
= d_A^{(12-r)/40}m^{-1}|I|\leq d_A^{-1/40}
\]
and by the same token
\[
\sum_{u\notin  U_2}\|(t_j)_{i\in I_{u}'}\|_2^2
\leq d_A^{-1/40}.
\]
Hence
\begin{align}
\notag
\sum_{u\notin U_1}\sum_{i\in I_u'',j\in I_u'}a_{i,j}\ve_{i,j}s_it_j
&\leq \|(a_{i,j}\ve_{i,j})_{i,j}\|\sum_{u\notin U_1}\|(s_i)_{i\in I_{u}''}\|_2\|(t_j)_{j\in I_{u}'}\|_2
\\
\notag
&\leq \|(a_{i,j}\ve_{i,j})_{i,j}\|
\Bigl(\sum_{u\notin U_1}\|(s_i)_{i\in I_{u}''}\|_2^2\Bigr)^{1/2}
\Bigl(\sum_{u\leq l}\|(t_j)_{j\in I_{u}'}\|_2^2\Bigr)^{1/2}
\\
\label{eq:estYkmA2}
&\leq d_A^{-1/80}\|(a_{i,j}\ve_{i,j})_{i,j}\|
\end{align}
and
\begin{align}
\notag
\sum_{u\in U_1\setminus U_2}\sum_{i\in I_u'',j\in I_u'}a_{i,j}\ve_{i,j}s_it_j
&\leq \|(a_{i,j}\ve_{i,j})_{i,j}\|
\sum_{u\notin U_2}\|(s_i)_{i\in I_{u}''}\|_2\|(t_j)_{j\in I_{u}'}\|_2
\\
\notag
&\leq \|(a_{i,j}\ve_{i,j})_{i,j}\|
\Bigl(\sum_{u\leq l}\|(s_i)_{i\in I_{u}''}\|_2^2\Bigr)^{1/2}
\Bigl(\sum_{u\notin U_2}\|(t_j)_{j\in I_{u}'}\|_2^2\Bigr)^{1/2}
\\
\label{eq:estYkmA3}
&\leq d_A^{-1/80}\|(a_{i,j}\ve_{i,j})_{i,j}\|.
\end{align}

Bounds \eqref{eq:estYkmA0}-\eqref{eq:estYkmA3} yield
\[
\Ex\max_{m}Y_{m,r}\leq \Ex\max_{k}Z_k+2d_A^{-1/80}\Ex\|(a_{i,j}\ve_{i,j})_{i,j}\|
\]
and the assertion follows by Lemma \ref{lem:maxconnYA} and Proposition \ref{prop:normgraphA1}.
\end{proof}

\medskip

\noindent
\textbf{Acknowledgements.}
Part of this work was carried out while the author was visiting the Hausdorff Research 
Institute for Mathematics, Univerity of Bonn. The hospitality of HIM and of the organizers of the program \textit{Synergies between modern probability, geometric analysis and stochastic
geometry} is gratefully acknowledged.
The author would like to thank Marta Strzelecka for careful readings of various versions of 
the manuscript and numerous  valuable comments, Ramon van Handel for pointing out the reference 
\cite{BiL} and an anonymous referee for many thoughtful comments and for suggesting most of the examples presented in the paper.



\vspace{0.4cm}

\begin{thebibliography}{99}
\bibitem{AL} R.~Adamczak and R.~Latała, 
\textit{Tail and moment estimates for chaoses generated by symmetric random variables with logarithmically concave tails,}
Ann. Inst. Henri Poincaré Probab. Stat. \textbf{48} (2012), 1103--1136.

\bibitem{ACKM} N.~Agarwal,  K.~Chandrasekaran, A.~Kolla and V.~Madan,
\textit{On the expansion of group-based lifts,}
SIAM J. Discrete Math. \textbf{33} (2019), 1338--1373.


\bibitem{AGZbook} G.~W.~Anderson, A.~Guionnet, and O.~Zeitouni,
\textit{An introduction to random matrices,}
Cambridge Studies in Advanced Mathematics, vol. 118, Cambridge University Press,
Cambridge, 2010.

\bibitem{BBvH} A.~S.~Bandeira, M.~T.~Boedihardjo, and R.~van Handel, 
\textit{Matrix concentration inequalities and free probability}, 
Invent. Math. \textbf{234} (2023), 419--487. 

\bibitem{BvH} A.~S.~Bandeira and R.~van Handel, 
\textit{Sharp nonasymptotic bounds on the norm of random matrices with independent entries,} 
Ann. Probab. \textbf{44} (2016), 2479--2506.

\bibitem{BeL} W.~Bednorz and R.~Latała,
\textit{On the boundedness of Bernoulli processes},
Ann. of Math. (2) \textbf{180} (2014), 1167--1203.

\bibitem{BiL} Y.~Bilu and N.~Linial,
\textit{Lifts, discrepancy and nearly optimal spectral gap},
Combinatorica \textbf{26} (2006), 495--519.

\bibitem{Bo} C.~Bordenave
\textit{ A new proof of Friedman's second eigenvalue theorem and its extension to random lifts},
Ann. Sci. Éc. Norm. Supér. (4) \textbf{53} (2020), 1393--1439.


\bibitem{DMS} S.~J.~Dilworth and S.~J.~Montgomery-Smith, 
\textit{The distribution of vector-valued Rademacher series}, 
Ann. Probab. \textbf{21} (1993), 2046--2052.

\bibitem{Hi} P.~Hitczenko,
\textit{Domination inequality for martingale transforms of a Rademacher sequence}, 
Israel J. Math. \textbf{84} (1993), 161--178.

\bibitem{LS} R.~Latała and W.~Świątkowski, 
\textit{Norms of randomized circulant matrices,} 
Electron. J. Probab. \textbf{27} (2022), Paper No. 80, 23pp.

\bibitem{LvHY} R.~Latała, R.~van Handel, and P.~Youssef, 
\textit{The dimension-free structure of nonhomogeneous random matrices}, 
Invent. Math. \textbf{214} (2018), 1031--1080.

\bibitem{Ledoux-Talagrand} M. Ledoux and M. Talagrand, 
\textit{Probability in Banach spaces. Isoperimetry and processes,}
Ergebnisse der Mathematik und ihrer Grenzgebiete (3), vol. 23, Springer-Verlag, Berlin, 1991.

\bibitem{MOP} S.~Mohanty, R.~O'Donnell and P.~Paredes,
\textit{Explicit near-Ramanujan graphs of every degree},
SIAM J. Comput. 51 (2022), no. 3, STOC20-1--STOC20-23.

\bibitem{MS} S.~J.~Montgomery-Smith,  
\textit{The distribution of Rademacher sums},
Proc. Amer. Math. Soc. \textbf{109} (1990), 517--522.

\bibitem{Seginer} Y.~Seginer, 
\textit{The expected norm of random matrices}, 
Combin. Probab. Comput. \textbf{9} (2000), 149--166.

\bibitem{vH} R.~van Handel, 
\textit{On the spectral norm of Gaussian random matrices,} 
Trans. Amer. Math. Soc. \textbf{369} (2017), 8161--8178.

\bibitem{vH2017_survey} R.~van Handel, 
\textit{Structured random matrices}, Convexity and concentration, IMA Vol. Math.
Appl., vol. 161, Springer, New York, 2017, pp. 107--156.

\end{thebibliography}
\end{document}